\documentclass[15pt]{article}

\usepackage{amsmath,amssymb}

\Large
\def\supp{\operatorname{supp}}

\def\e{\varepsilon}

\newtheorem{exam}{Example}[section]
\newtheorem{thm}{Theorem}[section]
\newtheorem{df}{Definition}[section]

\newtheorem{lem}{Lemma}[section]

\newtheorem{rem}{Remark}[section]
\numberwithin{equation}{section}

\newcommand{\qed}{\hfill$\Box$}

\newcommand{\dstl}{\displaystyle}

\newcommand{\wt}[1]{\widetilde{#1}}

\newcommand{\wbs}[1]{\,\overline{\!{#1}}}

\newcommand{\wbml}[1]{\;\!\overline{\:\!\!{#1}\;\!\!}}

\newcommand{\ch}[1]{\chi\kern-.05em\lower1ex\hbox{$\scriptstyle{#1}$}}
\newcommand{\chb}[1]
{\wbml{\chi}\kern-.02em\lower1ex\hbox{$\scriptstyle{#1}$}\;\!}

\newcommand{\vLambda}{{\mit\Lambda}}

\newcommand{\vOmega}{{\mit\Omega}}

\newcommand{\bs}{\:\!\!\setminus\:\!\!}

\newcommand{\Lap}{{\mit{\Delta}}}

\newcommand{\R}{{\bf{R}}}

\newcommand{\Rn}{{\bf R}\kern-0.08em\lower-0.75ex\hbox{$\scriptscriptstyle n$}}
\newcommand{\Rd}[1]{{\bf R}\kern-0.1em\lower-0.75ex\hbox{$\scriptscriptstyle{#1}$}}
\newcommand{\Sd}{S\kern0.0em\lower-0.75ex\hbox{$\scriptscriptstyle d\;\!\!-\:\!\!1$}}
\newcommand{\Sn}{S\kern0.0em\lower-0.75ex\hbox{$\scriptscriptstyle n\;\!\!-\:\!\!1$}}
\newcommand{\Cn}{{\bf C}\kern-0.01em\lower-0.75ex\hbox{$\scriptscriptstyle n$}}
\newcommand{\el}[1]{{\ell}\kern0.02em\lower-0.75ex\hbox{$\scriptscriptstyle{#1}$}}
\newcommand{\EL}[1]{{L}\kern0.0em\lower-0.75ex\hbox{$\scriptscriptstyle{#1}$}}
\newcommand{\EA}[1]{{A}\kern0.0em\lower-0.75ex\hbox{$\scriptscriptstyle{#1}$}}
\newcommand{\Tst}{T\kern0.0em\lower-0.58ex\hbox{$\scriptscriptstyle *$}}

\newcommand{\wtM}[1]{{\:\wt{\;\!\!\!M\;\!\!}}
\kern0.1em\lower-0.9ex\hbox{$\scriptstyle{#1}$}}
\newcommand{\wbM}[1]{{\;\!\wbs{\:\!\!M\:\!\!}}
\kern0.14em\lower-0.9ex\hbox{$\scriptstyle{#1}$}}

\title{ On general  Caffarelli-Kohn-Nirenberg type inequalities\\ involving non-doubling
weights \\ in the case of $p=1$  }
\author{   Toshio Horiuchi}

\begin{document}
\maketitle


\begin{abstract}

 We study the  Caffarelli-Kohn-Nirenberg type inequalities in the case of $p=1$  and 
generalize them  
adopting   weight functions
 $w(|x|)$ on $\R^n$ with  $w(t)\in  {W}(\R_+)$.  Here
${W}(\R_+)$ is  a  general class of  weight functions on $\R_+$ including  non-doubling weights like $e^{\pm 1/t}$. 
\footnote{2010{\it Mathematics Subject Classification.} Primary 35J70; Secondary 35J60.
\\ 
{\it Key words and phrases}. the  Caffarelli-Kohn-Nirenberg type inequalities, non-doubling weights, the  case of $p=1$ 
\\
 This research was partially supported by Grant-in-Aid for Scientific Research (No. 21K03304, No. 26K06858).}

\end{abstract}

\section{Introduction }
The  main purpose of the present paper is 
to study the Caffarelli-Kohn-Nirenberg type inequalities, which are abbreviated as  the CKN-type inequalities.
The CKN-type inequalities were introduced in \cite{CKN} as  multiplicative interpolation inequalities, but here we refer to the simple weighted Sobolev inequalities. There is a great deal of research in that case alone and we also studied in  \cite{ho1,ho2, hk3,CH,CH2}. Recently in \cite{ho4, ah4,ho3}, we revisited the CKN-type inequalities
in the case that  $p>1$ and   established   the CKN-type inequalities
involving  non-doubling weights.  
Furthermore, the CKN-type  inequalities, which differ greatly  in critical and non-critical cases,  were successfully  unified by adopting a new framework.
In the present paper we will proceed to study   the CKN-type inequalities in the case of $p=1$. \par
Let  $\R_+= (0,\infty)$.
First we  define a class of weight functions $W(\R_+)$ which  is a slight modification of the space introduced in \cite{ho4} 
to suit our purpose
(c.f.  Remark \ref{remark1.2}). 
By $C^{0,1}_{\rm loc}(\R_+)$ we  denote  the space  of  all locally Lipschitz continuous functions on $\R_+$.
Here, local Lipschitz continuousness of $v(t)$ means that $v(t)$ is Lipschitz continuous over each compact set of $\R_+$. 
\begin{df}\label{D0}
 \begin{enumerate}\item
 For $a\in [0,\infty]$ we  define
 \begin{equation}{W}_a(\mathbf R_+)=
\{ w\in C^{0,1}_{\rm loc}({\R}_+):   w>0, \lim_{t\to+0}w(t)  = a\, \}\end{equation}
 and  
 \begin{equation}{W}(\R_+)=\cup_{a\in [0,\infty]} {W}_a(\mathbf R_+).
\end{equation}
  
\item In particular we  set
  \begin{align}
&V(\R_+)= {W}_0(\mathbf R_+)\cup {W}_\infty(\mathbf R_+).
\end{align}

\end{enumerate}
\end{df}
We will  study  the CKN-type inequalities for  weights $w(t)\in V(\R_+)$ mainly in the subsequent. 
If $w(t)\in {W}(\R_+)\setminus V(\R_+)$,  then  $w(t)\in C^{0,1}_{\rm loc}({\R}_+)\cap C([0,\infty))$ and $0<\dstl\lim_{t\to+0}w(t)<\infty$,  hence
 $w(t)$ behaves tamely as $t\to +0$. 
 
\par
We introduce a kind of monotone rearrangement of  weight functions $w(t)\in V(\R_+)$.
\begin{df}\label{D2} 
\begin{enumerate}\item For $w(t)\in {W}_0(\R_+)$ and $0<\eta\le \infty$, we define
\begin{equation}
\varphi_w(t;\eta) = \begin{cases}
\dstl\inf_{ t\le s\le \eta} w(s)\quad &(0\le t\le\eta),\\
w(\eta) &(\eta\le t).
\end{cases}
\end{equation}
Here $\varphi_w \in {W}_0(\R_+)$ and $\varphi_w$   is  called  the largest  increasing function  with respect to $w$ 
statisfying
$\varphi_w(t;\eta)\le w(t) $ for $0\le t\le \eta.$
\item
For $w(t)\in {W}_\infty(\R_+)$ and $0<\eta\le \infty$, we define
\begin{equation}
\psi_w(t;\eta) = \begin{cases}\dstl \inf_{0\le s\le t}w(s)\qquad &(0\le t\le \eta),\\
 \dstl \inf_{0\le s\le \eta}w(s)\qquad &(\eta\le t).
\end{cases}
\end{equation}
Here   $\psi_w\in{W}_\infty(\R_+)$ and  $\psi_w$ is called the largest  decreasing function with respect to $w$ satisfying  $\psi_w(t)\le w(t) $ for $0\le t\le \eta.$
\item
For $w(t)\in V(\R_+)$ and  $0<\eta\le \infty$,  we define
\begin{equation} v_w(t;\eta)= \begin{cases}  \varphi_w(t;\eta),\qquad &w(t)\in {W}_0(\R_+),\\ \psi_w(t;\eta),
\qquad& w(t)\in {W}_\infty(\R_+).
\end{cases}
\end{equation}
\item
For  $w(t)\in V(\R_+)$,  $q\ge 1$  and $0<\eta\le \infty$, we  define $V_w^q(t;\eta)\in L^\infty_{\rm loc}((0,\eta])$ by
\begin{equation} V^q_w(t;\eta)= \begin{cases} \frac{d}{dt}( \varphi_w(t;\eta))^q,\qquad &w(t)\in {W}_0(\R_+),\\-\frac{d}{dt}( \psi_w(t;\eta))^q,
\qquad &w(t)\in {W}_\infty(\R_+).
\end{cases}
\end{equation}
\item

$\varphi_w(t;\eta)$, $\psi_w(t;\eta)$, $v_w(t;\eta)$ and $ V_w^q(t;\eta)$  are  abbreviated as $\varphi_w(t)$,  $\psi_w(t)$, $v_w(t)$ and $ V_w^q(t)$
respectively.

\end{enumerate}
\end{df}
\begin{rem}\label{remark1.1}
Since $\varphi_w(t)$ and $\psi_w(t)$ are locally Lipschitz continuous on $\R_+$ for $w(t)\in V(\R_+)$, 
they 
 are differentiable  a.e. on $\R_+$. In particular  
  the first order derivatives  in  the  distribution sense 
  $\varphi'_w(t), \psi'_w(t): \R_+
\to\R{\,}$  coincide with those   in  the  classical sense a.e. on 
$ \R_+$, i.e. 
$$\varphi'_w(t)=D\varphi_w(t){\,}\Big({\:\!\!}
 =\lim_{h\;\!\to\;\!0}\dfrac{1}{t}(\varphi_w(t+h)-\varphi_w(t)\:\!\!)\Big){\,\,\,}
 \mbox{ for a.e. }t\in \R_+.
$$
In particular we  see that $\varphi_w'(t), \psi_w'(t)\in L^\infty_{\rm loc}((0,\eta))$.

\end{rem}

Then we  study   inequalities  of  the following type: For $1\le q<\infty$, $0<\eta\le\infty$  and $w\in V(\R_+)$,  
there exists a positive  number $ C=C(q,\eta, w)\ge 1$ such  that  
 we  have   for any $ u(x)\in C_c^\infty(B_\eta\setminus \{0\}) $
\begin{equation}
\int_{B_\eta} |\nabla u(x)| w(|x|)|x|^{1-n}\,dx\ge  C \left(  \int_{B_\eta} 
 |u(x)|^q V^q_w(|x|)\, |x|^{1-n}\,dx \right)^{1/q},     \label{1.4'}
\end{equation}
where  $B_\infty=\R^n$ and 
$B_\eta= \{x\in \R^n: |x|<\eta \}$.
It is worth saying that  the  inequality (\ref{1.4'})  does  not  hold  unconditionally, unless $n=1$ or $q=1$. 
In order to study the validity of (\ref{1.4'}) for each $w\in V(\R_+)$,
we  introduce  {\bf the  non-degenerate condition }(\ref{K-condition}) in  Section 3, which  controls the behavior of $w$ near $t=0$, 
and then we  make  clear   the validity of  (\ref{1.4'})  under  (\ref{K-condition}) as  Theorem \ref{thm4.2}.
Roughly speaking,  (\ref{K-condition})  assures that  $w(t)$ does  not behave so  badly as $t\to +0$, and hence the function  
 $K(r)$  given by Definition (\ref{K(r)})  is bounded away from $0$. On  the contrary  if $\dstl\lim_{r\to +0}K(r)=0$ is  assumed, 
 then  by Theorem \ref{Thm4.3}  the  inequality (\ref{1.4'})  is not valid. In  Theorem \ref{Thm4.4} we also characterize a  set of  weight functions 
for  which  (\ref{K-condition}) is violated. 
The proofs of Theorem \ref{thm4.2}, Theorem \ref{Thm4.3} and Theorem \ref{Thm4.4} are given in Section{3}, Section 4 and Section 5 respectively.

\section{Main results}
First  we state  a  result in the one dimensional case which is rather simple.
\begin{thm}\label{Thm4.1} Assume  that
$1\le q<\infty$, $0<\eta\le \infty$ and 
 $w(t)\in V(\R_+)$.
Then, 
there exists a positive  number $C=C(q, \eta, w)\ge 1$ such  that  
 for any $u\in C_c^1((0,\eta))$
we  have
\begin{equation}
\int_0^\eta |u'(t)|w(t)\,dt\ge C \left(\int_0^\eta
 |u(t)|^q V^q_w(t)\,dt\right)^{1/q}. \label{4.1Thm4.1}
\end{equation}
Moreover 
either if  $w(t)\in {W}_0(\R_+)$ or if $w(t)\in {W}_\infty(\R_+)$, $\eta=\infty$ and $ \dstl\lim_{t\to\infty}w(t)=0$,
 then the best constant $C$ in (\ref{4.1Thm4.1}) equals $1$.
\end{thm}
\begin{rem}
If $w(t)\in V(\R_+)$ is  monotone, then  (\ref{4.1Thm4.1}) simply becomes 
\begin{equation}
\int_0^\eta |u'(t)|w(t)\,dt\ge C \left(\int_0^\eta
 |u(t)|^q \left|(w(t)^q)'\right|\,dt\right)^{1/q},   \qquad u\in C_c^1((0,\eta)) . \label{ineq.prop1.1}
\end{equation}
\end{rem}

Theorem \ref{Thm4.1} can be derived from classical one-dimensional inequalities with some modifications, but for the sake of self-completion we give a  direct proof in \S3 (c.f.\cite{m}).
In order to  proceed to the $n$-dimensional case, let  us prepare  more notations.
\begin{df}\label{dfa1}
For $v\in   C^{0,1}_{\rm loc}({\R}_+)$
we define the followings:
\begin{align}&
 Z[v]=\{{\;\!}r\in\R_+ :
 \mbox{$v$ is differentiable at $r$ and $ v'(r)\ne 0$}{\;\!}\},
\\ 
&
 Z_{0}^{}[v]=\{{\;\!}r\in\R_+ :
 \mbox{$v$ is differentiable at $r$ and ${\;\!}v'(r)=0$}{\;\!}\}.
\end{align}
\end{df}

\begin{df}\label{K(r)} 
For $w(r)\in V(\R_+)$ we set
\begin{equation}
\begin{split}K(r)
= \left|\frac{w(r) }{ r w'(r)} \right|\qquad  ( r\in Z[v_w]),  
\end{split}
\end{equation}
where $v_w$ is defined by Definition \ref{D2},3.
\end{df}
 Now  we introduce {\bf the non-degenerate condition} (NDC) on $K(r)$ which assures  that
 $K(r)$ is  bounded away from $0$ as $r\to + 0$. 
\begin{df}[\bf the non-degenerate condition]  Let $\eta>0$ and $w\in V(\R_+)$. A weight function $w$ is said to satisfy   the non-degenerate condition (NDC) if    
\begin{equation}C_0:=\inf_{ r\in (0,\eta] \cap  Z[v_w]} K(r)>0.\tag{NDC}
\label{K-condition} \end{equation}
\end{df}
\medskip
Then we state the $n$-dimensional CKN-type inequality  as a natural extension of 
 Theorem \ref{Thm4.1} and
 the classical CKN-type inequalities in  Appendix, and  the proof  will be  given in \S 4.
\par\medskip 
\begin{thm}\label{thm4.2} Let $n>1$, $1\le q<\infty$,  $ 0 \le 1 -1/q \le1/n$ and $\eta>0$. 
Assume  that $w(r)\in V(\R_+)$.
Moreover assume  that if $1<q$, $K(r)$ satisfies  (\ref{K-condition}).
 Then, there exists a positive number $C=C(q,\eta,w)$ such  that  we have  for any $u\in C_c^\infty(B_\eta\setminus \{0\})$ 
 \begin{equation}
\int_{B_\eta} |\nabla u(x)| w(|x|)|x|^{1-n}\,dx\ge  C \left(  \int_{B_\eta} 
 |u(x)|^q V^q_w(|x|)\,|x|^{1-n}\,dx \right)^{1/q}.\label{6.3}
\end{equation}

 \end{thm}
\par\medskip 
Conversely  we  have  the following which is proved in \S 5:
\begin{thm}\label{Thm4.3}
 Let $n>1 $, $1< q <\infty$, $ 0 \le 1-1/q \le1/n$ and $\eta>0$. 
 Assume  that $w(r)\in V(\R_+)$.
If  
\begin{equation}\dstl \lim_{\e\to+0}\sup_{ r\in (0,\e] \cap  Z[v_w]}K(r)=0,\label{limsup} 
\end{equation}
then the inequality (\ref{6.3}) does  not  hold.
\end{thm}

Roughly speaking,  either if
  $w$    vanishes   infinitely at  the origin or if
 $w$    blows up   infinitely at  the origin, then   (NDC) is violated.
 To explain more accurately, we introduce the following notion.
\par\medskip
\begin{df} For  $w(r)\in V(\R_+)$ we define the following:
\begin{enumerate}\item
For $w(r)\in W_0(\R_+)$, $w(r)$  is  said to  vanish   in infinite order at  the origin,  if and only if 
for some $C>0$ and 
 for an arbitrary positive integer $m$ there exists  positive number  $r_m$  satisfying $r_m\to 0 $ as $m\to\infty$ such  that  we have
\begin{equation} w(r_m)\le C(r_m)^m. \label{vanishinfinite}
\end{equation}
\item
For $w(r)\in W_\infty(\R_+)$,
$w(r)$ is  said to blow up   in infinite order at  the origin,  if and only if
 for some $C>0$ and 
 for an arbitrary positive integer $m$ there exists a positive number  $r_m$  satisfying $r_m\to 0 $ as $m\to\infty$ such  that  we have
\begin{equation} w(r_m)\ge  C(r_m)^{-m}. \label{blowinfinite}
\end{equation}
\end{enumerate}
\end{df}

By  Definition \ref{D2} and Remark \ref{remark1.1} we immediately have  the following:
\begin{lem}\label{lemma2.2}  \begin{enumerate}
\item 
Assume  that $w(r)\in {W}_0(R_+)$. Then $w(t) $ vanishes  in infinite order at  the origin,  if and only if
$\varphi_w(r) $  vanishes  in infinite order at  the origin.
\item
Assume  that $w(r)\in {W}_\infty(R_+)$. Then $w(t) $ blows up  in infinite order at  the origin,  if and only if
$\varphi_w(r) $ blows up  in infinite order at  the origin.
\end{enumerate}
\end{lem}
\par\noindent{\bf Proof:} 
Assume that $w(r) $ vanishes  in infinite order at  the origin. Since $\varphi_w(r)\le w(r)$  for $t\in [0,\eta]$, so $\varphi_w(r)$ does. Conversely assume that $\varphi_w(r) $ vanishes  in infinite order at  the origin, namely,  we have
(\ref{vanishinfinite}) for some sequence of positive numbers $\{r_m\}$.
 Since $\varphi_w(r)$ is a constant on each component $Z_0[\varphi_w]$, we can assume $\varphi_w(r_m)= w(r_m),\; m=1,2,\cdots$. Hence the assertion 1 follows.
The  assertion 2 can be shown in a  similar way.
\qed 
\par\medskip\noindent
Then we have  the following which is proved in \S6:
\par
\begin{thm} \label{Thm4.4} Let $w(r) \in V(\R_+)$.
If
  $w$    vanishes  in  infinite order at  the origin or if $w$   blows 
 up    in  infinite order at  the origin, then     (NDC) is violated.
\end{thm}
To make the theorem easier to understand, we present  typical examples.
\begin{exam} Let $\alpha>0$. 
When either $w(r)=e^{-r^{-\alpha}} \in {W}_0(\R_+)$ or $w(r)=e^{r^{-\alpha}}\in {W}_\infty(\R_+)$, \par $K(r)= O(r^\alpha)$ as $r\to +0$.
\end{exam}
\begin{exam}
Let  $1\le q <\infty$ and $0<\eta$.
 Let $w(r)= r^{\gamma}$. 
\begin{enumerate}\item
If $\gamma> 0$, then $w(t)\in {W}_0(\R_+)$ and we  have 
$
 K(r)=1/\gamma.$
\item
If $\gamma<0 $, then $w(r)\in {W}_\infty(\R_+)$ and we  have 
$ K(r)=-1/\gamma.$ 
\end{enumerate}
\end{exam}

\begin{rem}\label{remark1.2}\begin{enumerate}\item
In \cite{ho4}, 
we  have established the CKN-type inequalities for  $p>1$ with $C^1$-weight functions, and they
 remain valid for weight functions in $W(\R_+)$ defined by (1.1).  
 \item
According to \cite{ho4}, if  we define
  subclasses ${P}(\R_+)$ and ${Q}(\R_+)$ by
 \begin{equation}  \begin{cases}&{P}(\mathbf R_+)= \{ w(t)\in {W}(\R_+) : \,  w(t)^{-1} \notin L^1((0,\eta)) \, \text{ for some} \, \eta >0\},\\
 &Q(\mathbf R_+) =\{ w(t)\in { W}(\R_+) :  \, w(t)^{-1}\in L^1((0,\eta)) \, \text{ for any } \, \eta >0 \},
 \end{cases}\label{counterpart}
 \end{equation}
then, we see that
  ${W}(\R_+)= {P}(\R_+)\cup{ Q}(\R_+).$ 
Here  we  have  the relations
    $${P}(\R_+)\subset {W}_0(\R_+), \; {W}_{\infty}(\R_+)\subset W(\R_+)\setminus {W}_0(\R_+) \subset {Q}(\R_+)\; \text { and }\;
{W}_0(\R_+)\cap {Q}(\R_+)\neq \phi.$$ 
\end{enumerate}
\end{rem}
\begin{rem}Interestingly
the inequality
 (\ref{4.1Thm4.1}) can  be  proved 
by taking the limit $(p\to 1+0)$  in the one-dimensional CKN-type inequalities for $p>1$   in  \cite{ho4}(Theorem 3.1 with $w^{p-1}=v^p \; 
(v\in V(\R_+)))$.
 \end{rem}

\section{Proof of Theorem \ref{Thm4.1}}
Without loss of generality we assume  that $u(t)\ge 0$.
First we consider the case that $w(t) \in {W}_0(\R_+)$. 
Recall $ V^q_w(t)= \frac{d}{dt}( \varphi_w(t))^q\ge 0$ a.e. in $(0,\eta].$
For $u(t)\in C_c^1((0,\eta))$ we set $f(t)= |u'(t)|$. Noting that  $0\le u(t) \le \int_t^\eta f(s)\,ds$, we have
\begin{align*}& \left(\int_0^\eta|u(t)|^q V^q_w(t)\,dt\right)^{1/q}\le 
\left( \int^\eta _0 \left|  \int_t ^\eta f(s) \,ds\right|^q V^q_w(t)\,dt\right)^{1/q}\\
&=\left( \int^\eta _0 \left|  \int_0^\eta f(s) V^q_w(t)^{1/q} \chi_{[t,\eta]}(s)\,ds\right|^q \,dt\right)^{1/q}
\le  \int^\eta _0 \left(  \int_0^\eta |f(s)|^q V^q_w(t) \chi_{[t,\eta]}(s)\,dt\right)^{1/q} \,ds\\
&=\int_0^\eta |f(s)| \left(  \int_0^s V^q_w(t)\,dt\right)^{1/q}\,ds
\le  \int_0^\eta |f(s)| (\varphi_w(s)^q-\varphi_w(0)^q)^{1/q}\,ds\\
&\le  \int_0^\eta |u'(s)| \varphi_w(s)\,ds\le \int_0^\eta |u'(s)| w(s)\,ds.
\end{align*}
Therefore we have (\ref{4.1Thm4.1}) with $C\ge 1$.\par
Secondly  we assume that $w\in {W}_\infty(\R_+)$. Recall $ V^q_w(t)= -\frac{d}{dt}( \psi_w(t))^q\ge 0\; (0<t\le \eta)$.  
Then   by noting that
$ u(t)\le \int_0^t f(s)\,ds $ with $f(t)=|u'(t)|$, we have
\begin{align*}&\left(\int_0^\eta
 |u(t)|^q V^q_w(t)\,dt\right)^{1/q}\le
\left( \int^\eta _0 \left|  \int_0 ^t f(s) \,ds\right|^q V^q_w(t)\,dt\right)^{1/q}\\
&=\left( \int^\eta _0 \left|  \int_0^\eta f(s) V^q_w(t)^{1/q} \chi_{[0,t]}(s)\,ds\right|^q \,dt\right)^{1/q}
\le  \int^\eta _0 \left(  \int_0^\eta |f(s)|^q V^q_w(t)\chi_{[0,t]}(s)\,dt\right)^{1/q} \,ds\\
&=\int_0^\eta |f(s)| \left(  \int_s^\eta V^q_w(t) \,dt\right)^{1/q}\,ds
\le  \int_0^\eta |f(s)| (\psi_w(s)^q-\psi_w(\eta)^q)^{1/q}\,ds\\
&\le  \int_0^\eta |f(s)| \psi_w(s)\,ds\le  \int_0^\eta |f(s)| w(s)\,ds.
\end{align*}
Thus we  have (\ref{4.1Thm4.1}) with $C\ge 1$. 
\par\medskip \noindent
Proof of Optimality: 
First we  assume that  $w(t)\in {W}_0(\R_+)$,  then we  have  (\ref{4.1Thm4.1}). 
By the previous argument, we  also have   for  an arbitrary  measurable function $f(t)\ge0$ 
\begin{align*} 
C\left( \int^\eta _0 \left|  \int_t ^\eta f(s) \,ds\right|^q V^q_w(t)\,dt\right)^{1/q}
\le  \int_0^\eta f(s) v_w(s)\,ds.\end{align*}
Now  we assume  that 
 $\supp f\subset (x, x+h)$ with $x\in (0,\eta)\cap (Z[v_w]\cup Z_0[v_w]), 0<h< \eta-x$.
Then 

\begin{align*}
\int_x^{x+h} f(t) v_w(t)\,dt &\ge C
\left( \int_{0} ^x \left(  \int_x^{x+h} f(s) \,ds\right)^q V^q_w(t)\,dt\right)^{1/q}
=  C  \varphi_w(x)  \int_x^{x+h} f(s) \,ds.
\end{align*}
Setting  $ f(t)= \varphi_w(t)^{-1}$ in $(x,x+h)$, we have
$$  1\ge  C  \frac {\varphi_w(x)}{h}\int_x^{x+h} \varphi_w(s)^{-1} \,ds.
$$
By letting $h\to 0$ we  see $C\le1$.
This proves the assertion. \par
Secondly  we  assume that $\eta=\infty$ and   $\dstl\lim_{t\to \infty}w(t)=0$. 
We   have   for  an arbitrary  measurable function $f(t)\ge0$  having a compact support 
\begin{align*} 
C\left( \int^\infty_0 \left|  \int_0 ^t f(s) \,ds\right|^q V^q_w(t)\,dt\right)^{1/q}
\le  \int_0^\infty f(s) v_w(s)\,ds.\end{align*}
Now  we assume  that 
 $\supp f\subset (x-h, x)$ with $x\in  Z[v_w]\cup Z_0[v_w]$ and $0<h<x$.
Then 
\begin{align*}
\int_{x-h}^{x} f(t) v_w(t)\,dt &\ge C
\left( \int_{x} ^\infty\left(  \int_{x-h}^{x} f(s) \,ds\right)^q V^q_w(t)\,dt\right)^{1/q}
=  C  \psi_w(x)  \int_{x-h}^{x} f(s) \,ds.
\end{align*}

Setting  $ f(t)= \psi_w(t)^{-1}$ in $(x-h,x)$, we have
$$  1\ge  C  \frac {\psi_w(x)}{h}\int_{x-h}^{x} \psi_w(s)^{-1} \,ds.
$$
By letting $h\to 0$ we  see $C\le1$.
This proves the assertion. 
\qed

\section{Proof of Theorem \ref{thm4.2}}
Assume  that  $w(t)\in V(\R_+)$,
$0<\eta<\infty$ and  $ 0 \le 1 -1/q \le1/n$.
By $\mu_1$ we  denote ($1$-dimensional) Lebesgue measure. Then we prepare the following.
\begin{lem}\label{Lemma2.1}
Assume that $w(t)\in V(\R_+)$. 
\begin{enumerate}\item
We have for an arbitrary compact set $K\subset \R_+$
$$ \mu_1(\R_+ \setminus (Z(v_w)\cup Z_0(v_w))=0\quad \mbox{  and }\quad \mu_1( v_w( Z_0(v_w))\cap K)=0.$$
\item We  have $v_w(r)=w(r)$ for all $r\in Z[v_w]$.  In particular  we have $v'_w(r)=w'(r)$ for all $r\in Z[v_w]$.
\item
Let  $\tilde{\eta}= v_w(\eta)$. Then,
$\rho=v_w(r)$ is  invertible on $(0,\eta]\setminus Z_0[v_w]$ and  
the inverse mapping  $v_w^{-1}: \rho\in (0,\tilde{\eta}]\setminus v_w(Z_0[v_w]) \mapsto r\in (0,\eta]\setminus Z_0[v_w]$ is differentiable a.e. to obtain $(v_w^{-1})'(\rho)= 1/ v_w'(r)$. 
\end{enumerate}

\end{lem}
\noindent{\bf Proof: } $1.$ 
From Remark \ref{remark1.1} 
for $w(r)\in V(\R_+)$ both $\varphi_w(r)$  and $\psi_w(r)$ are monotone continuous piecewise $C^1$ function on $\R_+$. Hence $ v_w(r)$ is a locally Lipschitz continuous function on $\R_+$ and 
 is differentiable  a.e. on $ \R_+$. 
In particular we have  ${\;\!}{{\mu}_{1}^{}
(\R_+{\;\!\!}\bs{\;\!\!}({ Z}[v_w]\cup{ Z}_{0}^{}[v_w])\:\!\!)}
=0{\;\!}$.  Moreover 
 we have for  an arbitrary 
 Borel set  ${\:\!}A\subset{(0,{\;\!\!}\infty)}{\:\!}$ 
$$\dstl{
 {\mu}_{1}^{}(v_w(A)\cap K\:\!\!) \le {\int_{A\cap(v_w)^{-1}(K) }}|v_w'{\;\!\!}(r)|{\:\!}dr
}.$$
Here, ${\;\!}{v_w(A)}={\{{\;\!}v_w(r) :  r\in A{\;\!}\}}{\;\!}$ and $(v_w)^{-1}(K)=\{r :  v_w(r)\in K\}$. Particularly we have
$$\dstl{
 {\mu}_{1}^{}(v_w({ Z}_{0}^{}[v_w])\cap K\:\!\!)
 \le {\int_{{ Z}_{\;\!\!0}^{}[v_w] \cap(v_w)^{-1}(K)}}|v_w'{\;\!\!}(r)|{\:\!}dr
 ={\int_{{ Z}_{\;\!\!0}^{}[v_w]\cap(v_w)^{-1}(K)}}|Dv_w(r)|{\:\!}dr=0
}.$$
\par\noindent $2.$  This follows direct  from Definition \ref{D2}.
\par\noindent $3.$  $v_w(r)$ is monotone  and $v'_w(r)\neq 0$ over $(0,\eta] \setminus Z_0[v_w]$, hence the assertion follows.
\par
\hfill \qed

We  use a polar coordinate system $x= r\omega $ for  $r= |x|$ and $\omega\in S^{n-1}$.
By $\Lap_{\Sn}^{}$  we  denote  the Laplace-Beltrami operator on a unit sphere $\dstl{
 S^{\:\!n\:\!-1}}$, and  by $dS$ we  denote    surface elements on  $\dstl{
 S^{\:\!n\:\!-1}}$. Then  a gradient operator $\vLambda$ on 
 $\dstl{
 S^{\:\!n\:\!-1}
}$ is defined by  
\begin{equation}\dstl{
 {\int_{\;\!\!\Sn}}
 (-{\;\!}\Lap_{\Sn}^{}\xi_1){\:\!}\xi_2{\;\!}dS
 ={\int_{\;\!\!\Sn}}\vLambda{\:\!}\xi_1{\cdot}\vLambda{\:\!}\xi_2{\;\!}dS
 {\,\,\,}\mbox{ for }\xi_1,\xi_2\in C^{\:\!2}(S^{\:\!n\:\!-1})
}.
\end{equation}
Here we note  that 
\begin{equation}\dstl{
 \Lap{\:\!}\xi_1=\dfrac{1}{r^{\:\!n-1}\;\!\!}
{\partial_r}\big(r^{\:\!n-1}{\partial_r{\:\!}\xi_1}\big)
 +{\;\!\!}\dfrac{1}{r^{2}\;\!\!}\Lap_{\Sn}^{}\xi_1,{\,\,}
 |\nabla \xi_1|^{2}=\big|{\partial_r\xi_1}\big|^{2}{\;\!\!}
 +{\;\!\!}\dfrac{1}{r^{2}\;\!\!}|\vLambda{\:\!}\xi_1|^{2}
},
\end{equation}
where
$
{\partial_r\xi_1}(x)={x}/{|x|}{\cdot}\nabla \xi_1(x).
$
By   a polar coordinate system, the inequality  (\ref{6.3}) is  transformed to the  following: 
There exists a positive number $C=C(q,\eta,w)$ such  that  we have for any  $u(x)= u(r\omega)\in C_c^\infty(B_\eta\setminus \{0\})$ 
\begin{equation}
\begin{split}
\int_{S^{n-1}} \,dS \int _0^\eta\left( (\partial_r u)^2+ \frac{(\Lambda u)^2}{r^2}\right)^{1/2} w(r)\,dr\ge  C \left( \int_{S^{n-1}} \,dS\int _0^\eta
|u|^q V^q_w(r)\,dr \right)^{1/q}.
\end{split}\label{5.6}
\end{equation}

Define a change of  variables  in harmony with subclasses  as follows:
\begin{enumerate}\item
If   $w(r)\in {W}_0(\R_+)$,  then we  set $\, v_w(r)= \rho,$ $\;\,( 0<r\le \eta)$.
\item
If  $w(r)\in {W}_\infty(\R_+)$, then we  set
$\, v_w(r)=1/\rho,$ $\;\,( 0<r\le \eta)$.
\end{enumerate}

First  we consider the case 1 to show (\ref{5.6}). 
From Lemma \ref{Lemma2.1}, $v_w(r)$ is invertible on $(0,\eta]\setminus Z_0[v_w]$.
We  employ a polar coordinate system $x= r\omega $ for  $r= |x|$ and $\omega\in S^{n-1}$. 
By  the change  of variable $v_w(r)=\rho$ for $r\in (0,\eta]\setminus Z_0[v_w]$ , we  have for $u(r\omega)\in C^\infty_c(B_{\eta}\setminus\{0\})$
\begin{align*}
\int _0^\eta
|u(r\omega)|^q V^q_w(r)\,dr= \int_{(0,\eta]\setminus Z_0[v_w]}|u(r\omega)|^q \frac{d}{dr}(v_w(r))^q\,dr
=\int_{(0,\tilde\eta]\setminus v_w(Z_0[v_w])}|U(\rho\omega)|^q d(\rho^q)
\end{align*}
where $U(\rho\omega)= u(v_w^{-1}(\rho)\omega) \in C(B_{\tilde \eta}\setminus\{0\})$ for $ \rho\in  (0,\tilde\eta]\setminus v_w(Z_0[v_w])$ and $v_w(\eta)=\tilde \eta$. \par
Since  $U(\rho\omega)$ is differentiable a.e. on $\{ (\rho,\omega): \rho\in  (0,\tilde\eta]\setminus v_w(Z_0[v_w]),\; \omega \in  S^{n-1}\}$,  we have 
\begin{align*}
\int _0^\eta \left( (\partial_r u)^2+ \frac{(\Lambda u)^2}{r^2}\right)^{1/2} w(r)\,dr &
\ge \int _{(0,\eta])\setminus Z_0[v_w]} \left( (\partial_r u)^2+ \frac{(\Lambda u)^2}{r^2}\right)^{1/2} v_w(r)\,dr\\
=& 
\int _{(0,\tilde\eta])\setminus v_w(Z_0[v_w])} \left( (\partial_\rho U)^2+ K(r)^2\frac{(\Lambda U)^2}{\rho^2}\right)^{1/2} \rho\,d\rho,
\end{align*}
where $r=v_w^{-1}(\rho)$ for $ \rho\in  (0,\tilde\eta]\setminus v_w(Z_0[v_w])$  and $K(r)$ is given by
\begin{equation}
\begin{split}K(r)
=\left|\frac{v_w(r) }{ r v_w'(r)} \right|= \left|\frac{w(r) }{ r w'(r)} \right|\quad  ( r\in (0,\eta]\cap Z(v_w)).  \end{split}\label{4.4}
\end{equation}
 Since   $\mu_1((0,\eta] \setminus (Z(v_w)\cup Z_0(v_w))=0$  and $\mu_1((0,\tilde\eta]\cap v_w(Z_0[v_w]))=0$ hold, together with a density argument w.r.t. $U(y)$ in $C^\infty_c(B_{\tilde\eta}\setminus \{0\})$,
the inequality (\ref{5.6}) is  reduced  to  the following: For  $U(y)=U(\rho \omega) \in C^\infty_c(B_{\tilde\eta}\setminus \{0\})$
\begin{equation}
\begin{split}
\int_{S^{n-1}} \,dS &\int _0^{\tilde \eta}\left( (\partial_\rho U)^2+H(\rho)^2 \frac{(\vLambda U)^2}{\rho^2}\right)^{1/2} \rho\,d\rho \ge  C \left( \int_{S^{n-1}} \,dS\int _0^{\tilde \eta}
 |U|^q   \,d(\rho^q)\right)^{1/q}, 
\end{split}\label{5.7}
\end{equation}
where $\tilde\eta=v_w(\eta)$ and 
\begin{equation}
H(\rho)= K(v_w^{-1}(\rho))=\left|\frac{\rho \;(v_w^{-1})'(\rho)}{v_w^{-1}(\rho)}\right|
 \;\mbox { for } \;  \rho\in (0,\tilde{\eta}]\cap v_w(Z[v_w])). \label{4.6}
 \end{equation}
If $q=1$ holds, then (\ref{5.7}) follows direct from the classical Hardy inequality, hence  we assume  that $q>1$. 
From (\ref{K-condition}),  the left hand side of (\ref{5.7}) is estimated from below in the following way. For $\rho=|y|$,

 \begin{align*}
\int _{B_{\tilde \eta}} \left( (\partial_\rho U)^2+H(\rho)^2 \frac{(\vLambda U)^2}{\rho^2}\right)^{1/2}\rho^{2-n}\,dy
&\ge \min(C_0,1)
\int _{B_{\tilde \eta}}\left( (\partial_\rho U)^2+\frac{(\vLambda U)^2}{\rho^2}\right)^{1/2} \rho^{2-n}\,dy\\ &= \min(C_0,1)
\int _{B_{\tilde \eta}}|\nabla_y U|\rho^{2-n}\,dy
\\
&\ge \min(C_0,1)S^{1,q;1} \left( \int_{B_{\tilde \eta}} 
 |U|^q   \rho^{q -n} \,dy\right)^{1/q}. \label{}
\end{align*}
In the last step we used the CKN type inequality (\ref{1.4}) with $\gamma=1$. This  proves   (\ref{5.7}) with $C=\min(C_0,1)S^{1,q;1}  q^{-1/q}$.
\par\medskip
Secondly  we  consider the case 2.  
By  the change  of variable $v_w(r)=1/\rho$ for $r\in (0,\eta]\setminus Z_0[v_w],$  we  have for $u(r\omega)\in C^\infty_c(B_{\eta}\setminus\{0\})$
\begin{align*}
\int _0^\eta
|u(r\omega)|^q V^q_w(r)\,dr= -\int_{(0,\eta]\setminus Z_0[v_w]}|u(r\omega)|^q \frac{d}{dr}(v_w(r))^q\,dr
=\int_{(0,\tilde\eta]\setminus \widetilde {v_w}(Z_0[v_w])}|U(\rho\omega)|^q d(-\rho^{-q})
\end{align*}
where $\widetilde {v_w}=1/{v_w}$, $U(\rho\omega)= u(v_w^{-1}(1/\rho)\omega) \in C(B_{\tilde \eta}\setminus\{0\})$ for $ \rho\in  (0,\tilde\eta]\setminus \widetilde {v_w}(Z_0[v_w])$ and $\tilde \eta=\widetilde {v_w}(\eta)\; (=1/v_w(\eta))$. 
Since  $U(\rho\omega)$ is differentiable a.e. on $\{ (\rho,\omega): \rho\in  (0,\tilde\eta]\setminus \widetilde{v_w}(Z_0[v_w]),\; \omega \in  S^{n-1}\}$,  we have 
\begin{align*}
\int _0^\eta \left( (\partial_r u)^2+ \frac{(\Lambda u)^2}{r^2}\right)^{1/2} w(r)\,dr &
\ge \int _{(0,\eta])\setminus Z_0[v_w]} \left( (\partial_r u)^2+ \frac{(\Lambda u)^2}{r^2}\right)^{1/2} v_w(r)\,dr\\
=& 
\int _{(0,\tilde\eta])\setminus \widetilde{v_w}(Z_0[v_w])} \left( (\partial_\rho U)^2+ K(r)^2\frac{(\Lambda U)^2}{\rho^2}\right)^{1/2} \rho^{-1}\,d\rho,
\end{align*}
where $r=v_w^{-1}(1/\rho)$ for $ \rho\in  (0,\tilde\eta]\setminus \widetilde{v_w}(Z_0[v_w])$  and $K(r)$ is given by (\ref{4.4}).
 Since   $\mu_1((0,\eta] \setminus (Z(v_w)\cup Z_0(v_w))=0$  and $\mu_1((0,\tilde\eta]\cap \widetilde{v_w}(Z_0[v_w]))=0$ hold, together with a density argument w.r.t. $U(y)$ in $C^\infty_c(B_{\tilde\eta}\setminus \{0\})$,
the inequality (\ref{5.6}) is  reduced  to  the following: For  $U(y)=U(\rho \omega) \in C^\infty_c(B_{\tilde\eta}\setminus \{0\})$
\begin{equation}
\begin{split}
\int_{S^{n-1}} \,dS &\int _0^{\tilde \eta}\left( (\partial_\rho U)^2+H(\rho)^2 \frac{(\Lambda U)^2}{\rho^2}\right)^{1/2}\rho^{-1}\,d\rho 
\ge  C \left( \int_{S^{n-1}} \,dS\int _0^{\tilde \eta}
 |U|^q   \,d(-\rho^{-q})\right)^{1/q},
\end{split}\label{5.9}
\end{equation}
where 
\begin{equation}
H(\rho)= K( v_w^{-1}(1/\rho))= \left|\frac{(v_w^{-1})'(1/\rho)}{\rho\;v_w^{-1}(1/\rho)}\right|\; \mbox{ for }\;
\rho\in (0,\tilde{\eta}]\cap \widetilde{v_w}(Z[v_w])) \label{4.8}
\end{equation}
As in the previous case, it suffices to show (\ref{5.9}) for $U(y)=U(\rho\omega)\in C^{\infty}_c(B_{\tilde \eta}\setminus \{0\})$.
Again we  assume  that $q>1$ and 
we   have
 \begin{align*}
\int _{B_{\tilde \eta}} \left( (\partial_\rho U)^2+H(\rho)^2 \frac{(\vLambda U)^2}{\rho^2}\right)^{1/2} \rho^{-n}\,dy
&\ge  
\min(C_0,1)
\int _{B_{\tilde \eta}}\left( (\partial_\rho U)^2+\frac{(\vLambda U)^2}{\rho^2}\right)^{1/2} \rho^{-n}\,dy\\ &= \min(C_0,1)
\int _{B_{\tilde \eta}}|\nabla_y U|\rho^{-n}\,dy
\\
&\ge \min(C_0,1)S^{1,q; -1} \left( \int_{B_{\tilde \eta}} 
 |U|^q   \rho^{-q -n} \,dy\right)^{1/q}. \label{}
\end{align*}
In the last step we used the CKN type inequality (\ref{1.4}) with $\gamma=-1$. This  proves   (\ref{5.9}) with $C=\min(C_0,1)S^{1,q;1} q^{-1/q}$. 
We  note  that  $S^{1,q;1}=S^{1,q;-1}=S_{1,q}= \omega_n^{1-1/q}q^{1/q}$ holds by Theorem \ref{thm1.3}.
Here  by $\omega_n$ we  denote  a surface  area of an $n$-dimensional unit ball.
\qed
\begin{rem}
If $C_0\ge 1$, then  $C=S^{1,q;1}q^{-1/q}=S^{1,q;-1}q^{-1/q}=S_{1,q}q^{-1/q}= \omega_n^{1-1/q}$ is the best constant. 
In fact, by Theorem \ref{thm1.3} with $n\ge 2$  and $\gamma =\pm1$, 
$S^{1,q; 1}=S^{1,q;-1}=S^{1,q; 1}_{\rm rad}=S^{1,q; -1}_{\rm rad}=S_{1,q} =\omega_n^{1-1/q}q^{1/q}$ holds.
Then one can assume  $U\in C_c^\infty(B_{\tilde{\eta}}\setminus \{0\})_{\rm rad}$ so that we have $\vLambda U\equiv 0$. Therefore the assertion is now clear.
\end{rem}

\section{Proof of Theorem \ref{Thm4.3}}
\par\noindent{\bf Proof:} 
{\bf Case1.} 
First we assume that 
$w(r)\in {W}_0(\R_+)$. Then    $ v_w(r)=\varphi_w(r)$ and  $ \dstl \lim_{r\to+0}v_w(r)=0 $. 
By a change of  variable  $v_w(r)=\rho$ for $r\in (0,\eta]\setminus Z_0[v_w]$,  $\tilde\eta=v_w(\eta)$  and 
$\tilde\e=v_w(\e)$,
the  inequality (\ref{6.3}) is  equivalent to 
 (\ref{5.7}). 
From the assumption (\ref{limsup}) and (\ref{4.6}) we   have
\begin{equation}\dstl \lim_{\e\to+0}\sup_{ r\in (0,\e]\cap Z[v_w]}K(r)=
\dstl \lim_{\e\to+0}\sup_{ \rho\in (0,\tilde\e)]\cap v_w(Z[v_w])}H(\rho) =0.
\label{limsup2}
\end{equation} 
Let  $B(\omega)\in L^1(S^{n-1})$ with 
$B(\omega) \notin L^q(S^{n-1})\; (q>1)$. 
Let $B_j(\omega)$ be a mollification of $B$ such  that $B_j(\omega)\in C^\infty(S^{n-1})$,
$B_j \to B  $ in $L^1(S^{n-1})$ but
\begin{equation}
\int_{S^{n-1}}|B_j(\omega)|^q\,dS \to \infty \quad (j\to\infty).\label{infty}
\end{equation}
 Let $\{ \e_j\}$ be  a sequence of  numbers  such that $0<\e_j<1$,    $\e_j\to 0$ as $j\to\infty$ 
and 
\begin{equation}
H(\rho)\cdot\int _{S^{n-1}} |\Lambda B_j(\omega)| \,dS \le 1\quad ( \rho\in (0, \e_j\tilde\eta]\cap(v_w(Z[v_w]) \setminus v_w(Z_0[v_w])),  \, j=1,2,3,\ldots).\label{6.25}
\end{equation}
We take and fix  an $ A(\rho)\in  C^\infty_c((0,  \tilde \eta))\setminus\{0\}$ satisfying  
\begin{equation} \int_0^{\tilde \eta} |\partial_\rho A(\rho)|\rho
\,d\rho =1. \label{6.26} \end{equation}
Define 
\begin{equation}A_j(\rho)=  \e_j^{-1}A(\rho/{\e_j})\qquad (j=1,2,3,\ldots ). \label{Aj}\end{equation}
Then,  we see that for $j=1,2,3,\ldots$ 
\begin{equation}\begin{cases}& A_j(\rho)\in  C^\infty_c((0,\e_j \tilde \eta)),\\
 &\int_0^{\tilde \eta \e_j} |\partial_\rho A_j(\rho)|\rho
\,d\rho =\int_0^{\tilde \eta}  |\partial_\rho A(\rho)|\rho
\,d\rho =1,\\ 
&\int_0^{\tilde \eta\e_j}  |A_j(\rho)|^q\rho^{q-1}
\,d\rho=\int_0^{\tilde \eta} | A(\rho)|^q\rho^{q-1}
\,d\rho <\infty.
\end{cases}
\label{6.28} \end{equation} 
Then we define a sequence of test functions $U_j= A_j(\rho)\cdot B_j(\omega) \in  C^\infty_c((0,\tilde \eta))\times C^\infty(S^{n-1})$.
If we  show the following properties, then the assertion clearly follows:
\begin{align}\label{6.30}&
\int_{S^{n-1}} \,dS \int _0^{\tilde \eta\e_j}\left( (\partial_\rho U_j)^2+H(\rho)^2 \frac{(\Lambda U_j)^2}{\rho^2}\right)^{1/2} \rho\,d\rho <\infty,
\\
&\label{6.31}
 \left( \int_{S^{n-1}} \,dS\int _0^{\tilde \eta\e_j}
 |U_j|^q   \rho^{q-1} \,d\rho\right)^{1/q} \to \infty \quad \mbox{as }\quad j\to\infty.
\end{align}
From (\ref{infty}) and (\ref{6.28}) we have (\ref{6.31}),   hence
 it suffices to show (\ref{6.30}).
We note that
\begin{equation}\label{6.32}
 \int_0^{\tilde \eta\e_j} |\partial_\rho A_j(\rho)|\rho
\,d\rho \int_{S^{n-1} } |B_j(\omega)| \,dS= \int_{S^{n-1} } |B_j(\omega)| \,dS<\infty. \end{equation} 
From (\ref{6.25}), (\ref{6.25}) and  the fact $\mu_1( v_w(Z_0[v_w])\cap \supp A_j)=0$  we have 
\begin{equation}
\begin{split}
\int_0^{\tilde \eta\e_j} |A_j(\rho)| H(\rho)\,d\rho \int_{S^{n-1}} |\Lambda B_j(\omega)|\, dS
\le \int_0^{\tilde \eta\e_j} |A_j(\rho)| \,d\rho 
=\int_0^{\tilde \eta} | A(\rho)|
\,d\rho
<\infty
\end{split}\label{6.34}
\end{equation}
Since $(a^2+b^2)^{1/2}\le 2^{1/2}(a+b),\,(a,b\ge 0)$, we have (\ref{6.30}), hence the assertion is proved.
\par {\bf Case 2.}
Secondly we assume  that $w\in {W}_\infty(\R_+)$ and (\ref{limsup}).
From the assumption (\ref{limsup})  and (\ref{4.8}) we   have
\begin{equation}\dstl \lim_{\e\to+0}\sup_{ r\in (0,\e]\cap Z[v_w]}K(r)=
\dstl \lim_{\e\to+0}\sup_{ \rho\in (0,\tilde\e)]\cap \widetilde {v_w}(Z[v_w])}H(\rho) =0.
\label{limsup2}
\end{equation}

Let $B_j(\omega)\in C^\infty(S^{n-1}) \, (j=1,2,3,\ldots)$ be the same function as before. 
We take   an $ A(\rho)\in  C^\infty_c((0,  \tilde \eta))\setminus\{0\}$ satisfying  
\begin{equation} \int_0^{\tilde \eta} |\partial_\rho A(\rho)|\rho^{-1}
\,d\rho =1. \label{6.35} \end{equation}
Define 
\begin{equation}A_j(\rho)= \e_j A(\rho/{\e_j})\qquad (j=1,2,3,\ldots ). \label{Aj2}\end{equation}
Then,   we see that for $j=1,2,3,\ldots$
\begin{equation}\begin{cases}& A_j(\rho)\in  C^\infty_c((0,\e_j \tilde \eta)),\\
 &\int_0^{\tilde \eta \e_j} |\partial_\rho A_j(\rho)|\rho^{-1}
\,d\rho =\int_0^{\tilde \eta}  |\partial_\rho A(\rho)|\rho^{-1}
\,d\rho =1,\\ 
&\int_0^{\tilde \eta\e_j}  |A_j(\rho)|^q\rho^{-1-q}
\,d\rho=\int_0^{\tilde \eta} | A(\rho)|^q\rho^{-1-q}
\,d\rho <\infty.
\end{cases}
\label{6.37} \end{equation} 
Now we define a sequence of test functions $U_j= A_j(\rho)\cdot B_j(\omega) \in  C^\infty_c(0,\tilde \eta\e_j)\times C^\infty(S^{n-1})$.
If we can show the following properties, then the assertion  follows in a similar way:
\begin{align}\label{6.39}&
\int_{S^{n-1}} \,dS \int _0^{\tilde \eta\e_j}\left( (\partial_\rho U_j)^2+H(\rho)^2 \frac{(\Lambda U_j)^2}{\rho^2}\right)^{1/2} \rho^{-1}\,d\rho <\infty,
\\
&\label{6.40}
 \left( \int_{S^{n-1}} \,dS\int _0^{\tilde \eta\e_j}
 |U_j|^q   \rho^{-1-q} \,d\rho\right)^{1/q} \to \infty \quad \mbox{as }\quad j\to\infty.
\end{align}
Since (\ref{6.40}) follows direct from (\ref{infty}) and  (\ref{6.37}), it suffices to show (\ref{6.39}).
Again we note that
\begin{equation}\label{6.41}
 \int_0^{\tilde \eta\e_j} |\partial_\rho A_j(\rho)|\rho^{-1}
\,d\rho \int_{S^{n-1} } |B_j(\omega)| \,dS= \int_{S^{n-1} } |B_j(\omega)| \,dS<\infty. \end{equation} 
Then  we have
\begin{equation}
\begin{split}
\int_0^{\tilde \eta\e_j} |A_j(\rho)| & H(\rho)\rho^{-2}\,d\rho \int_{S^{n-1}} |\Lambda B_j(\omega)|\, dS\\
& \le \int_0^{\tilde \eta\e_j} | A_j(\rho)|\rho^{-2}
\,d\rho=\int_0^{\tilde \eta} | A(\rho)|\rho^{-2}
\,d\rho
<\infty\qquad  ( (\ref{6.25}) ).
\end{split}\label{6.43}
\end{equation}
Hence  we have (\ref{6.39}) as before, and the assertion is proved.
\qed

\section{Proof of Theorem \ref{Thm4.4}}
\par\noindent{\bf Proof:} 
First we treat the  case that $w(r)\in {W}_0(\R_+)$. Then  $v_w(r)=\varphi_w(r)$ and  $\dstl\lim_{r\to+0}v_w(r)=0$.
From Lemma \ref{lemma2.2},  $v_w$ vanishes in infinite order at  the origin.
Namely  we assume that for some $C >0 $ and for  an arbitrary positive number $m$ there exists a  positive $r_m$  such  that $r_m\to 0 $ as $m\to\infty$ and 
\begin{equation} v_w(r_m)\le C(r_m)^m. \label{6.46}
\end{equation}
Now we assume  on the contrary that for  some  positive numbers $C_0$ and (a small) $\eta$, 
\begin{equation}K(r) \ge C_0, \qquad 0\le r\le \eta.\label{assumption}\end{equation}
Then  $C_0 v_w'(r)/v_w(r) \le  1/r$ holds for $r\in (0,\eta]\setminus Z_0[v_w])$, 
hence  this  holds over $(0,\eta]$. By integrating the  both side over an interval $[r, \eta]$ we have 
\begin{equation} 
 v_w(\eta) \left(\frac{r}{\eta}\right)^{1/C_0} \le v_w(r), \qquad r \in (0,\eta]. \label{6.47}
\end{equation}
Then from (\ref{6.46}) we  have 
%
\begin{equation} v_w(\eta) \left(\frac{r_m}{\eta}\right)^{1/C_0}\le C(r_m)^m,  \qquad m=1,2,\cdots. \label{6.48}
\end{equation}
If $m$ is sufficiently large, then this does not hold, hence
 the assertion is proved by  a contradiction.
\par\medskip
Secondly  we treat the  case that $w(r)\in {W}_\infty(\R_+)$. 
Then  $v_w(r)=\psi_w(r)$ and  $\dstl\lim_{r\to+0}v_w(r)=\infty$.
From Lemma \ref{lemma2.2},  $v_w$ blows up in infinite order at  the origin.
%
Then we assume  that for some  $C>0$ and 
 for an arbitrary positive number $m$ there exists a positive $r_m$  such  that 
 $r_m\to 0 $ as $m\to\infty$ and
\begin{equation} w(r_m)\ge C(r_m)^{-m},  \qquad m=1,2,\cdots. \label{6.49}
\end{equation}
As  in the previous step  we  assume (\ref{assumption}).
Noting  that  $ K(r)= -v_w(r)/(r v_w'(r))$ for $r\in (0,\eta]\setminus Z_0[v_w]$, we have
$-C_0 v_w'(r)/v_w(r) \le  1/r$ holds for $r\in (0,\eta]\setminus Z_0[v_w]$, hence by integrating the  both side over an interval $[r, \eta]$ we have 
\begin{equation} 
 v_w(\eta) \left(\frac{r}{\eta}\right)^{-1/C_0} \ge v_w(r), \qquad r \in (0,\eta). \label{6.47}
\end{equation}
Then from (\ref{6.49}) we  have 
\begin{equation} v_w(\eta) \left(\frac{r_m}{\eta}\right)^{-1/C_0}\ge C(r_m)^{-m},  \qquad m=1,2,\cdots. \label{6.48}
\end{equation}
If $m$ is sufficiently large again, then this does not hold, hence
 the assertion is proved by  a contradiction.
{}\hfill\qed

\section{Appendix:   
The  non-critical CKN-type inequalities}
In the  non-critical case, the CKN-type  inequalities have the following form:
\begin{equation}\int_{\R^n}|\nabla u(x)| |x|^{1+\gamma -n}\, dx\ge
 S^{1,q;\gamma}\bigg(\int_{\R^n}|u(x)|^q |x|^{\gamma q-n}\, dx\bigg)^{1/q}, u \in 
C_{\rm c}^{\infty}(\R^n\setminus \{0\}), 
\label{1.4}
\end{equation}
where    $n \ge 1$  and $q,\gamma$ are real numbers satisfying 
 \begin{equation} \gamma \neq 0,\quad  q<\infty, \quad
 0 \le 1 -1/q \le1/n .\label{1.5}\end{equation} 
%

Here $S^{p,q;\gamma} $ is  called the  best constant and given by

\begin{align}
 S^{\:\!1\;\!\!,q\:\!;\:\!\gamma}
 &=\inf\{{\;\!}E^{\:\!1\;\!\!,q\:\!;\:\!\gamma}[u]\mid u\in
 C^{\:\!\infty}_{\rm c}(\R^{\;\!\!n}\bs{\;\!\!}\{0\}{\:\!\!}){\;\!\!}\bs{\;\!\!}\{0\}\}, \label{best}
\end{align}
where
\begin{equation}
 E^{\:\!1\;\!\!,q\:\!;\:\!\gamma}[u]
 =\dfrac{\int_{\R^n}|\nabla u(x)| |x|^{1+\gamma -n}\, dx}
 {\bigg(\int_{\R^n}|u(x)|^q |x|^{\gamma q-n}\, dx\bigg)^{1/q}}{\,\,\,}\mbox{ for }u\in
 C^{\:\!\infty}_{\rm c}(\R^{\;\!\!n}\bs{\;\!\!}\{0\}{\:\!\!}).
\end{equation}

 We  also define the  radial best constant as follows.
 \begin{df}
Let 
$\vOmega$ be  a radially symmetric domain. 
For any function space  $V(\vOmega)$ on 
$\vOmega$,   we set
\begin{equation}\dstl{
 V(\vOmega)_{\rm rad}^{}=\{{\;\!}u\in V(\vOmega)\mid\mbox{$u$ is radial}{\;\!}\}.
}\end{equation}
\end{df}
 Then we define
\begin{equation}
S^{\:\!1\;\!\!,q\:\!;\:\!\gamma}_{\rm rad}
 =\inf\{{\;\!}E^{\:\!1\;\!\!,q\:\!;\:\!\gamma}[u]\mid u
 \in C^{\:\!\infty}_{\rm c}(\R^{\;\!\!n}\bs{\;\!\!}\{0\}{\:\!\!})_{\rm rad}^{}{\!}
 \bs{\;\!\!}\{0\}\}. \label{symmetric}
\end{equation}
\begin{rem}\label{remark7.1}
Here we remark  that
the  best  constants $S^{1,q;\gamma} $ is invariant if the whole space $\R^n$ is  replaced by an  arbitrary bounded domain $\Omega$ containing the origin.  $ S^{\:\!1\;\!\!,q\:\!;\:\!\gamma}_{\rm rad} $ is  also invariant  if $\R^n$ is  replaced by a  radially symmetric domain $\Omega$.
For the detailed information see \cite{ho4}, \cite{hk3}.
\end{rem}

\begin{df}\label{df1.3}Let $\omega_n$ be a surface  area of an $n$-dimensional unit ball. For
$1\le q<\infty{\;\!}$,  we set
\begin{equation}
  S_{\;\!\!1\;\!\!,q}^{}
  = \omega_n^{1-1/q }q^{1/q}.\label{1.9}
\end{equation}
\end{df}

\begin{thm}[Symmetry]\label{thm1.3} Let $n\geq 1$. Assume  that
$1\le q <\infty$ and  $\tau_{\;\!\!1\;\!\!,q}^{}= 1-1/q\le1/n{\;\!}$. Then it holds that: 
\begin{enumerate}
\item
$\dstl{
 S^{\:\!1\;\!\!,q\:\!;\:\!\gamma}=S^{\:\!1\;\!\!,q\:\!;\:\!-\:\!\gamma},{\,\,}
 S^{\:\!1\;\!\!,q\:\!;\:\!\gamma}_{\rm rad}
 =S^{\:\!1\;\!\!,q\:\!;\:\!-\:\!\gamma}_{\rm rad} \hspace{1.58cm}
 \mbox{ for }\gamma\ne 0.
}$
\item
$\dstl{
 S^{\:\!1\;\!\!,q\:\!;\:\!\gamma}_{\rm rad}
 =S_{\;\!\!1\;\!\!,q}^{}|\gamma{\:\!}|^{1-\:\!\tau_{\;\!\!1\:\!\!,\;\!\!q}}
 \hspace{3.68cm} \mbox{ for }\gamma\ne 0.
}$
\item
$\dstl{
 S^{\:\!1\;\!\!,q\:\!;\:\!\gamma}=S^{\:\!1\;\!\!,q\:\!;\:\!\gamma}_{\rm rad}
 =S_{\;\!\!1\;\!\!,q}^{}|\gamma{\:\!}|^{1-\:\!\tau_{\;\!\!1\:\!\!,\;\!\!q}}
 \hspace{2.31cm} \mbox{ for }0<|\gamma{\:\!}|\le n-1, \, \, n>1.
}$
\end{enumerate}
\end{thm}

\par\bigskip
\noindent{{\bf Acknowledgments}}\par\medskip
The author would like to thank the referee for the precise advice in completing this paper.

\bigskip\noindent
{\large\bf Toshio Horiuchi\\Department of Mathematics\\Faculty of Science \\ Ibaraki University\\
Mito, Ibaraki, 310, Japan}\par\noindent
e-mail: toshio.horiuchi.math@vc.ibaraki.ac.jp
\par\noindent


\medskip


\begin{thebibliography}{BVxx}

%

\bibitem{ah4}  H. Ando, T. Horiuchi,  Generalized weighted Hardy's inequalities with compact perturbations, 
{\it Journal of Mathematical Inequalities,}  {\bf 18}, Number 1 (2024), pp 103-126.





\bibitem{CKN}  L. Caffarelli, R. Kohn,  L. Nirenberg, First order interpolation inequalities
with weights, {\it Compositio Math.},{\bf Vol. 53}, 1984, No. 3,
pp 259-275.  


\bibitem{CH}
N. Chiba, T. Horiuchi, 
On radial symmetry and its breaking in the Caffarelli-Kohn-Nirenberg type inequalities for $p=1$,
{\it Math. J. Ibaraki Univ.}, {\bf Vol.~47} (2015), pp 49--63.

\bibitem{CH2}
N. Chiba, T. Horiuchi, 
Radial symmetry and its breaking in the Caffarelli-Kohn-Nirenberg type inequalities for $p=1$,
{\it  Proc. Japan Acad., Ser. A, Math. Sci.}, {\bf Vol.~92} (2016), No.~4, pp 51--55. 

%
\bibitem{ho1} T. Horiuchi, The imbedding theorems for weighted Sobolev spaces, 
{\it Journal of
Mathematics of Kyoto University}, {\bf Vol. 29},  1989, pp 365-403.

\bibitem{ho2} T. Horiuchi, Best constant in weighted Sobolev inequality with weights being powers
of distance from the origin, {\it Journal of Inequality and Application}, {\bf Vol. 1},  1997, pp
275-292.


\bibitem{ho3} T. Horiuchi,
Hardy's inequalities with non-doubling weights and  sharp remainders,
 {\it SCMJ (in Editione Electronica) e-2022-2 Whole Number 35} (2022);
arXiv:2012.08766 [math.AP]

\bibitem{ho4} T. Horiuchi,
On general  Caffarelli-Kohn-Nirenberg type inequalities involving non-doubling
weights, {\it SCMJ (in Editione Electronica) e-2022-10 Whole Number 35} (2022)



\bibitem{hk3} T. Horiuchi, P. Kumlin, 
On the Caffarelli-Kohn-Nirenberg type inequalities involving Critical and Supercritical Weights, 
{\it Kyoto journal of Mathematics}, {\bf Vol. 52}, { No.4},  (2012),  pp 661-742.

%

\bibitem {m} V.G. Maz'ja, Sobolev spaces, {\it  
Springer}, 1985.


\end{thebibliography}
\end{document}